\documentclass[graybox]{svmult}

\usepackage{mathptmx}       
\usepackage{helvet}         
\usepackage{courier}        
\usepackage{type1cm}        
\usepackage{makeidx}         
\usepackage{graphicx}        
\usepackage{multicol}        
\usepackage[bottom]{footmisc}

\usepackage{amsmath}
\usepackage{amssymb}
\usepackage{psfrag} 
\usepackage{epsfig}

\newcommand{\ba}{\begin{array}}
\newcommand{\ea}{\end{array}}
\newcommand{\be}{\begin{equation}}
\newcommand{\ee}{\end{equation}}
\newcommand{\bd}{\begin{displaymath}}
\newcommand{\ed}{\end{displaymath}}
\newcommand{\bi}{\begin{itemize}}
\newcommand{\ei}{\end{itemize}}
\newcommand{\bn}{\begin{enumerate}}
\newcommand{\en}{\end{enumerate}}

\newcommand{\f}{\frac}

\newcommand*\samethanks[1][\value{footnote}]{\footnotemark[#1]}

\makeindex             

\begin{document}

\title*{On the scaling of entropy viscosity in high order methods}
\author{Adeline Kornelus\thanks{Supported in part
by NSF Grant DMS-1319054.
Any conclusions or recommendations expressed in this paper are those of the
author and do not necessarily reflect the views NSF.} and Daniel Appel\"{o}\samethanks}
\institute{Adeline Kornelus \at The Department of Mathematics and Statistics, University of New Mexico, Albuquerque, NM 87131, \email{kornelus@unm.edu}
\and Daniel Appel\"{o} \at The Department of Mathematics and Statistics, University of New Mexico, Albuquerque, NM 87131, \email{appelo@math.unm.edu}}
%
%
\maketitle


\abstract*{In this work, we outline the entropy viscosity method and discuss how the choice of scaling influences the size of viscosity for a simple shock problem. We present examples to illustrate the performance of the entropy viscosity method under two distinct scalings.}

\abstract{In this work, we outline the entropy viscosity method and discuss how the choice of scaling influences the size of viscosity for a simple shock problem. We present examples to illustrate the performance of the entropy viscosity method under two distinct scalings.}

\section{Introduction} \label{sec:intro}

Hyperbolic partial differential equations (PDE) are used to model various fluid flow problems. In the special case of 1-dimensional linear constant coefficient scalar hyperbolic problems, the solutions to these PDE are simply a translation of the initial data. However, for nonlinear problems the solution may deform, and as a result, shock waves can form even if the initial data is smooth \cite{leveque2002finite}.

In computational fluid dynamics, it is desirable that numerical methods capture shock waves and maintain a high accuracy for smooth waves. Low order methods have sufficient numerical dissipation to regularize shock waves but obtaining accurate solutions in smooth regions can become expensive. On the other hand, high order methods are capable of achieving high accuracy at a reasonable cost. Their low numerical dissipation enables such accuracy, but on the downside, it limits their ability to regularize shock waves. 

Various techniques have been implemented to capture shocks while maintaining high accuracy, at least away from shocks. There are two major classes of shock capturing techniques: shock detection techniques, where we find slope limiters \cite{leveque2002finite}, Essentially Non-Oscillatory (ENO) and Weighted ENO (WENO) \cite{shu1998essentially}, and artificial viscosity techniques, where we find filtering \cite{yee2007development, persson2006sub}, the PDE-based viscosity method \cite{johnson1990convergence}, the entropy viscosity method \cite{Guermond2011conservation}, among others.

In this work, we focus on the entropy viscosity method. In essence, the entropy viscosity method provides shock capturing without compromising the high accuracy away from the shock. An important advantage of this method is that it generalizes very easily to higher dimensions and unstructured grids.\\

As a model problem, we consider Burgers' equation
\be \label{burgers1d}
u_t + f(u)_x = 0,
\ee
where $f = \f{u^2}{2}$. Physically correct solutions to (\ref{burgers1d}) can be singled out by requiring that they satisfy an entropy inequality  such as
\be \label{ev_res}
r_{EV} = E_t + F_x \equiv \left( \f{u^2}{2}\right)_t + \left(\f{u^3}{3}\right)_x  \le 0. 
\ee
The entropy residual, $r_{EV}$, is zero wherever $u$ is smooth. If the solution $u$ contains a shock, then the entropy residual takes the form of a negative Dirac distribution centered at the location of the shock, $x_s$, i.e. $r_{EV} = -C \, \delta(x-x_s)$. The property that the entropy residual is unbounded at a shock was first used by Guermond and Pasquetti in \cite{Guermond2011conservation}, as a way to selectively introduce viscosity. The artificial viscosity, $\nu$, proposed in \cite{Guermond2011conservation}, defined as the minima of two viscosities   
\be \label{nutot}
\nu = \min (\nu_{\rm max},\nu_{EV}),
\ee 
becomes the coefficient of the viscous term in the viscous Burgers' equation, 
\be \label{burgers1d_visc}
u_t + f(u)_x = (\nu u_x)_x.
\ee

Here, $\nu_{\rm max}$ is the Lax-Friedrich viscosity whose size depends on discretization and the largest eigenvalue, $\lambda_{\rm LF}$, of the flux Jacobian, $\frac{Df(u)}{Du}$. The second viscosity $\nu_{EV}$ is proportional to the magnitude of the entropy residual (in fact, a discretization of the entropy residual) and will thus be zero (or small after discretization) away from discontinuities. In theory, the entropy residual becomes unbounded at a shock, numerically however, the entropy residual $r_{EV}$ remains bounded with the size of the residual depending on the discretization size. As we will see below, this subtle difference has consequences for how to choose the scaling of the viscosity terms in the entropy viscosity method.

On a grid with step size $h$, the second viscosity $\nu_{EV}$ can be expressed as
\be \label{nuev}
\nu_{EV}(x) = \alpha_{EV} h^\beta | r_{EV}(x) |,
\ee
with a parameter $\alpha_{EV}$ that requires tuning. In recent papers on entropy viscosity method, see e.g. \cite{guermond2011entropyhighorder,guermond2010entropy,guermond2011entropy,guermond2011suitable,zingan2013implementation}, the parameter $\beta$ is chosen to be 2, but the original paper \cite{guermond2008entropy} uses $\beta=1$. It is unclear to us why the later works prefer $\beta=2$. Here, we will present analysis and computational results that suggest the original scaling $\beta=1$ is a more natural choice. We note that the entropy residual is typically scaled by $\|E-\overline{E}\|_{\infty}$, with the over-bar indicating a spatial average, but as this quantity is roughly constant in the problems presented here, we omit it for brevity and reduced complexity.

The rest of the paper is organized as follows. In Section \ref{sec:num_method}, we describe different discretizations of (\ref{burgers1d_visc}) that we consider here, in Section \ref{sec:viscosity_analysis}, we present an analysis of how the entropy viscosity $\nu$ depends on the two viscosities, $\nu_{EV}$ and $\nu_{\rm max}$, under different scaling for a model problem. In Section \ref{sec:experiment}, we then conduct experiments with the entropy viscosity method where $\beta $ takes on values 1 or 2 and compare the results.

\section{Numerical methods}\label{sec:num_method}

We will consider the discretization of (\ref{burgers1d_visc}) by our conservative Hermite method \cite{KorAppJSC}, a standard discontinuous Galerkin (dG) method \cite{Hesthaven:2008fk} and a simple finite volume type discretization \cite{leveque2002finite}. For all the discretizations we let the domain $x_L \le x \le x_R $ be discretized by the regular grid $x_i = x_L + i\,h, \, i = 0,\ldots,n,$ $h = (x_R-x_L)/n$.

The degrees of freedom for the finite volume method is cell averages centered at the grid points. For the Hermite method, the degrees of freedom are the coefficients of node centered Taylor polynomials of degree $m$ and for the dG method, they are the ($m+1$) coefficients of element-wise (we take an element to be $\Omega_i = [x_{i-1},x_i]$) expansions in Legendre polynomials. For smooth solutions the spatial accuracy of the Hermite method is $2m+1$ and $m+1$ for the dG method. 

All three methods use the classic fourth order Runge-Kutta method to evolve the semi-discretizations in a method-of-lines fashion.

In the Hermite method, we evaluate the fluxes and their derivatives at the nodes (element edges) for the four stages in the RK method. Precisely, for the first stage we compute the slope $f^h_1 = \f{1}{2}\mathcal{T}[(u^h_1)^2]-\f{\nu}{h}\f{du^h_1}{dx}$ for the Taylor polynomial $u^h_1=u^h$ approximating the solution at the first stage. Here $\mathcal{T}[(u^h_1)^2]$ is the truncated polynomial multiplication of $u^h_1$ with itself and $\f{du^h_1}{dx}$ is the derivative of the polynomial. At the next stage, the solution is $u^h_2=u^h + \f{(\Delta t / 2)}{2} \f{df^h_1}{dx}$, the slope is $f^h_2 = \f{1}{2}\mathcal{T}[(u^h_2)^2]-\f{\nu}{h}\f{du^h_2}{dx}$  and so on. Once the stage slopes $f_s^h,\, s=1,\ldots,4$ and their spatial derivatives are known, we perform a Hermite interpolation to the element centers of the solution and the spatial derivatives of the stage slopes. These are then used to evolve the element centered Hermite interpolant of $u^h$ to $t = t_n + \Delta t /2$. As the Hermite interpolant is of higher degree than the original Taylor polynomial, we conclude a half-step by truncating it to the appropriate degree. To advance the solution a full time step, the half-step process is repeated starting from the element centers.

To handle the artificial viscosity in the dG method, we use the approach of Bassi and Rebay \cite{bassi1997high} with a Lax-Friedrichs flux for the advective term and alternating fluxes for the viscous term. The nonlinear terms are constructed explicitly and de-aliased by over-integration \cite{kirby2003aliasing}. 

For the finite volume method, we let $u_i \approx u(x_i)$ be a grid function approximating the solution and $f_{i+\f{1}{2}} = f_{i+\f{1}{2}}(u_i,u_{i+1})$ be an approximation to the flux at $x_{i+\f{1}{2}}$. To compute the time derivatives, we use the spatial approximation 
\be \label{fv_semidiscrete}
\f{du_i}{dt} \approx \f{f_{i+\f{1}{2}}- f_{i-\f{1}{2}}}{h},  
\ee
where 
\be
f_{i+\f{1}{2}}(u_i,u_{i+1}) = \f{1}{2} \left(\f{u_i+u_{i+1}}{2}\right)^2 - \left(\f{\nu_{i}+\nu_{i+1}}{2} \right) \frac{u_{i+1}-u_i}{h}. 
\ee
When $\nu_i = 0$, the above discretization is linearly stable (when paired with a suitable time-stepping method) but is not non-linearly stable, and we thus add artificial viscosity to stabilize it.

For all three discretizations, we approximate the time derivative of the entropy function, $E_t$, by a backward difference. This approach is explicit as we use the current solution to compute $E$ at the current time before evolving the solution in time. The residual (and hence the viscosity) is kept on each element / grid-point over each step. 

To approximate the entropy flux derivative $F_x$ using the Hermite method, we compute the derivative of the truncated polynomial multiplication $\mathcal{T}[u^h \mathcal{T}[(u^h)^2]]$ at the node. For the dG method, we evaluate the flux $F$ on a Legendre-Gauss-Lobatto (LGL) grid and differentiate it to get an approximation for $F_x$. The residual on an element is taken to be the maximum of the absolute value of the residual on the LGL grid. In the finite volume method $F_x$ is approximated by  
$$
\f{dF_i}{dx} = \f{F_{i+\f{1}{2}}- F_{i-\f{1}{2}}}{h}, \text{ where, } F_{i+\f{1}{2}} =  \f{1}{3} \left(\f{u_i+u_{i+1}}{2}\right)^3.
$$ 

We note that more sophisticated discretizations of the entropy residual could be considered. In particular, a higher order approximation to $r_{EV}$ would result in a higher rate of convergence for smooth solutions, but as we are mainly concerned with the scaling $\beta$, we did not pursue such discretizations here. In fact, in our experience, the results concerning the choice of scaling are not affected by the order of the accuracy of the approximation to $r_{EV}$. This will be discussed in Section \ref{sec:viscosity_analysis}.  

We also define $\nu_{\rm max}$ to be the classical Lax-Friedrich viscosity, which for Burgers' equation takes the form 
\be\label{nulf}
\nu_{\rm max} = \alpha_{\rm max} h \max | u |, 
\ee
where the maximum is taken globally. 

Finally, for the purpose of comparison we also present some results computed using the sub-cell resolution smoothness sensor of Persson and Peraire, \cite{persson2006sub}. The smoothness sensor compares the $L_2$ energy content of the highest (Fourier or expansion) mode with the total $L_2$ energy on an element and then maps its ratio (which is an indicator of the smoothness) into the size of the artificial viscosity. Precisely, if the approximate dG solution on an element is $u^h = \sum_{k = 0}^{m} \hat{u}_k P_k$, with $P_k$ being an orthogonal basis, we compute the smoothness as $s = \log_{10} (\|\hat{u}_m P_m \|^2 / \| u^h \|^2)$ and the viscosity as
\begin{equation*}
\nu = \left\{ \begin{array}{ll}
0 & s < s_0-\kappa, \\
\varepsilon_0 h & s > s_0+\kappa, \\
\f{\varepsilon_0 h}{2} \left(1 + \sin \left(  \f{\pi (s-s_0)}{2 \kappa}  \right) \right) & \text{otherwise}.
\end{array} \right.
\end{equation*}

When applied to the Hermite method, we first project the Taylor polynomials centered at two adjacent grid-points into an orthogonal Legendre expansion on the element defined by the grid-points and then proceed as above. 

\section{Impact of the $h$-scaling on the selection mechanism}
\label{sec:viscosity_analysis}

To study how the selection mechanism depends on the shock speed and the size of the jump, consider a solution of the Burgers' equation consisting of a Heaviside function $H$ with left state $u_l$ and right state $u_r$, given by 
\be \label{burger_shock_sol}
u(x,t) = u_l+\Delta u \  H\left(x-v_st\right).
\ee
This corresponds to a shock of size $|\Delta u|=|u_r-u_l|$ moving with speed $v_s=0.5(u_l+u_r)$. Solutions of the form (\ref{burger_shock_sol}) always has a negative $\Delta u$ value since Lax entropy condition for Burgers' equation dictates $u_l=f'(u_l)>v_s>f'(u_r)=u_r$.

For simplicity, we use the short hand notation $H$ for $H\left(x-v_st\right)$. A direct computation
\begin{align*}
u_t + \left(\f{u^2}{2}\right)_x & =  \left(-\f{u_l+u_r}{2}(\Delta u)H'\right)+\left((\Delta u)u_lH'+\f{(\Delta u)^2}{2}H'\right)\\
                               & = -\Delta u \left(\f{2u_l+\Delta u}{2}\right)H'+\Delta u\left(\f{2u_l+\Delta u}{2}\right)H'\\
                               & = 0,
\end{align*} 
shows that (\ref{burger_shock_sol}) is a solution of (\ref{burgers1d}). Further, it can be shown that the entropy residual (\ref{ev_res}) for (\ref{burger_shock_sol}) is    
\be
  r_{EV} = \f{(\Delta u)^3}{12}H'(x-x_s) = \f{(\Delta u)^3}{12}\delta(x-x_s).
\ee
That is, the size of the entropy residual grows with the cube of $\Delta u$. 

Now, by the properties that define the Dirac delta function $\delta$, we have 
\be
\int_{-\infty}^{\infty} \delta(x) dx = 1.
\ee
Thus, a consistent discretization of the Dirac delta function $\delta_0,...,\delta_n$ on a grid $x_0,...,x_n$ must obey the condition 
\be
\sum_{j=0}^{n-1}h_j \delta_j = 1,
\ee
where $h_j = x_{j+1}-x_j$. For any approximation with a finite width stencil, we must have $\delta_j\sim h^{-1}_j$ and we thus expect $r_{EV}$ to behave like $ (\Delta u)^3/h$ on a uniform grid. We therefore proceed with the analysis using the discrete approximation $r_{EV}=(\Delta u)^3/h$. Using this expression for $r_{EV}$, we estimate the viscosity $\nu$ by the minimum of   
\be
\nu_{EV} = {\alpha}_{EV} h^{\beta-1} |(\Delta u)^3|  \text{ and } \nu_{\rm max}= \alpha_{\rm max}h\max(|u_l|,|u_r|). 
\ee

The comparison between the size of $\nu_{EV}$ and $\nu_{\rm max}$ in various scenarios is reported in Table \ref{table:muE_mumax_comp}. If $\beta=2$, then the two viscosities $\nu_{EV}$ and $\nu_{\rm max}$ scale as $h$. For a problem with multiple shocks, the homogeneity in $h$-scaling introduces an additional difficulty in determining $ {\alpha}_{EV}$. Should it be chosen based on the largest or smallest shock? What if new shocks appear during the course of the computation? To avoid answering these questions, we instead consider $\beta=1$. Now $\nu_{EV} = \mathcal{O}(1)$ while $\nu_{\rm max} = \mathcal{O}(h)$, and the particular choice of $ {\alpha}_{EV}$ is thus irrelevant since as $h\to 0 $, the selection mechanism will eventually select $\nu_{\rm max}$ at the shocks. We will provide an example to illustrate the two-shock dilemma in Section \ref{sec:two_shocks}.

\begin{table}
\caption{Size of $\mu_E$ and $\mu_{\rm max}$ for different size of shock speed ($v_s$) with respect to the size of the jump ($\Delta u$) in the entropy viscosity method.}
\label{table:muE_mumax_comp}       
%
%
\begin{tabular}{p{2.4cm}p{3cm}p{5cm}}
\hline\noalign{\smallskip}
Case    & $\nu_{EV}$   & $\nu_{\rm max}$    \\
\noalign{\smallskip}\svhline\noalign{\smallskip}
{\bf $|v_S| \ll |\Delta u|$}     & $  {\alpha}_{EV}h^{\beta-1} |\Delta u|^3$ & $\alpha_{\rm max}h|v_s|$  \\     
{\bf $|v_s| \approx |\Delta u|$}  & $ {\alpha}_{EV}h^{\beta-1} |\Delta u|^3$ & $2\alpha_{\rm max}h|v_s|$\\
{\bf $|v_s| \gg |\Delta u|$}     & $ {\alpha}_{EV}h^{\beta-1} |\Delta u|^3$ & $0.5\alpha_{\rm max}h|\Delta u|$  \\         
\noalign{\smallskip}\hline\noalign{\smallskip}
\end{tabular}
\end{table}


\section{Experiments}\label{sec:experiment}

In this section, we describe the experiments and present a convergence study in $L_2$ norms, and also study the effects of the scaling in the entropy viscosity method on the convergence under grid refinement. For all the examples we solve Burgers' equation and vary the initial data. In each problem, we report the $L_2$-errors (the $L_1$-errors behaves quantitively similar). 

The solutions are obtained using the following methods: H1 and H2 refer to Hermite-entropy viscosity method for $\beta=1$ and $\beta=2$ respectively, DG1 and DG2 refer to dG-entropy viscosity method for $\beta=1$ and $\beta=2$ respectively, FV1 and FV2 refer to finite volume-entropy viscosity method with $\beta=1$ and $\beta=2$ respectively, DGP and HP refer to dG and Hermite method with smoothness sensor respectively.

The size of the time step is chosen close to the stability limit, which in the cases considered here results in the error being dominated by the spatial discretization.

\begin{figure}[htb]
  \begin{center} 
  \psfrag{XXXXXXXX}[][][0.8][0]{$h$} 
  \psfrag{YYYYYYYY}[][][0.8][0]{$L_2$-error} 
  \psfrag{AAA1}[][][0.6][0]{H1} 
  \psfrag{AAA2}[][][0.6][0]{H2} 
  \psfrag{AAA3}[][][0.6][0]{HP} 
  \psfrag{AAA4}[][][0.6][0]{DG1} 
  \psfrag{AAA5}[][][0.6][0]{DG2} 
  \psfrag{AAA6}[][][0.6][0]{DGP} 
  \psfrag{AAA7}[][][0.6][0]{FV1} 
  \psfrag{AAA8}[][][0.6][0]{FV2}
  \includegraphics[width=0.495\textwidth]{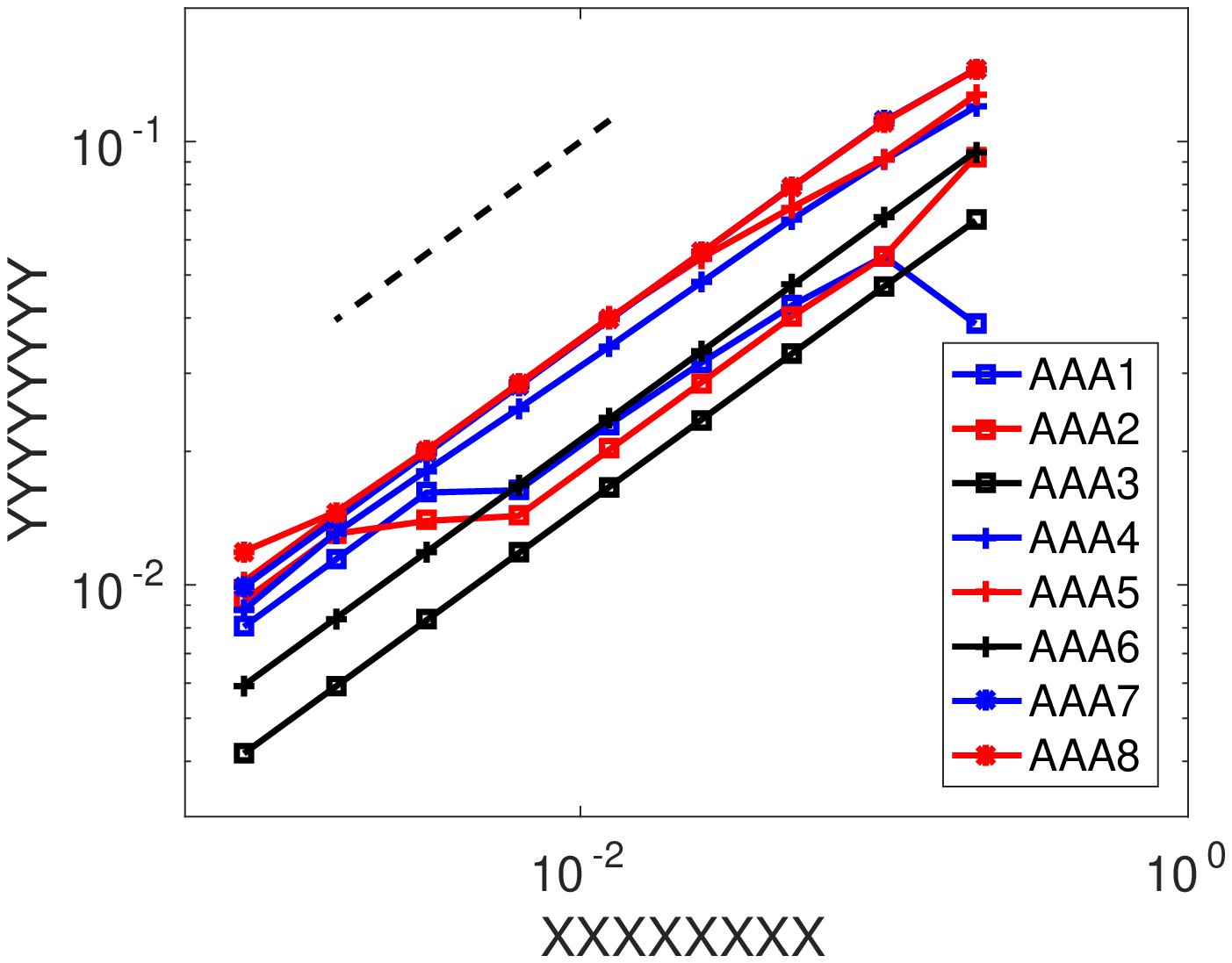}
  \psfrag{YYYYYYYY}[][][0.8][0]{} 
  \includegraphics[width=0.495\textwidth]{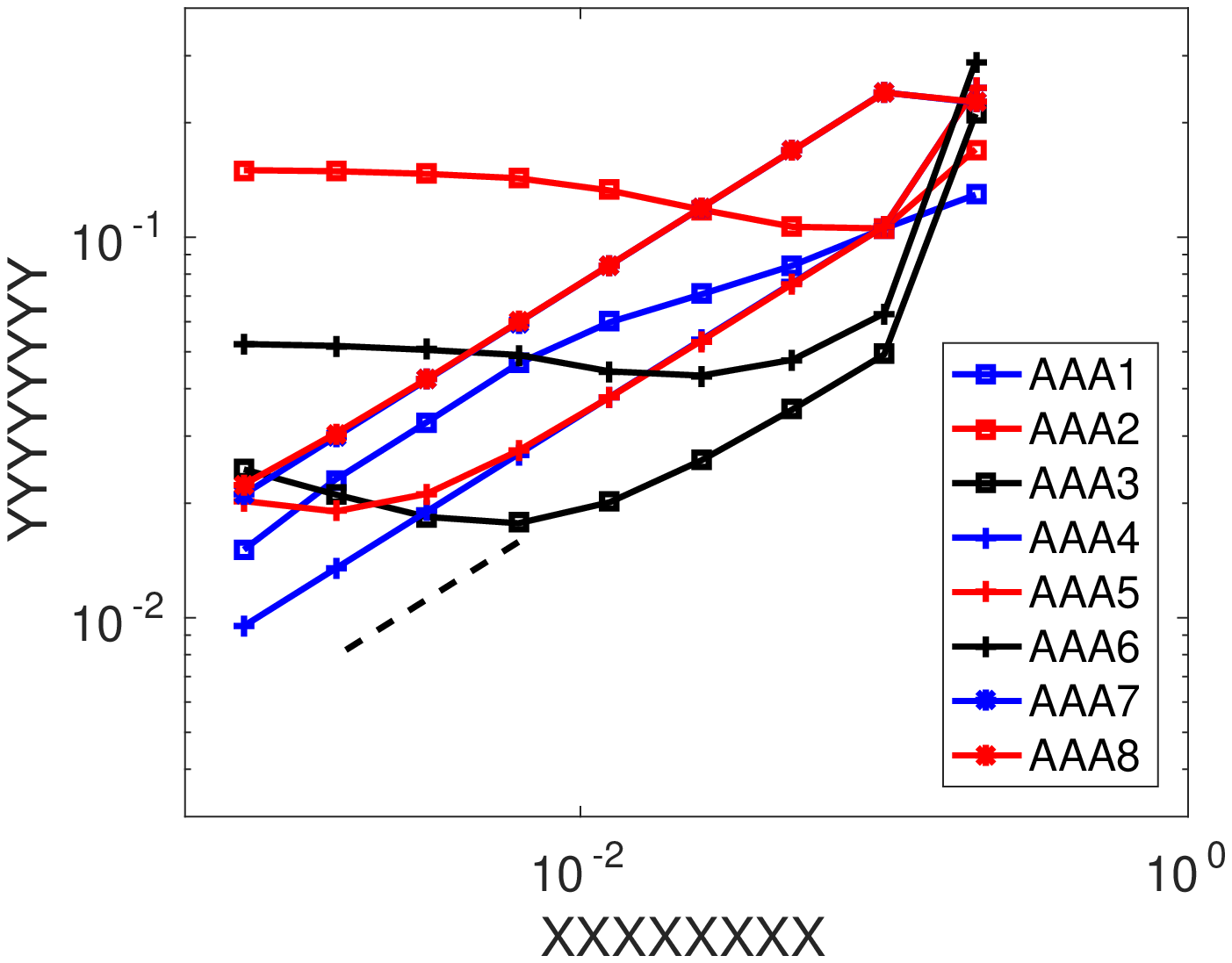}\hspace{0.1cm}
  \caption{Convergence of the different methods for stationary (left) and moving (right) shocks. \label{fig:APA}}
  \end{center}
\end{figure}

\subsection{A single shock}
\label{sec:single_shock}
In this example, we compute the solution to (\ref{burgers1d}) on the domain $D=[-1,1]$ with the initial data imposed as the exact solution 
\be
u(x,t) =
\begin{cases}
 -0.5+v_s,          & x \in [-1,v_st),\\
  \ \ \,  0.5+v_s,   & x \in [v_st,1],
\end{cases}
\ee
at time $t=0$. Here $v_s$ is the shock speed which we choose to be either $v_s=0$ corresponding to a stationary shock or $v_s=0.1$ corresponding to a moving shock.

We solve until time $t=1$ for the two different shock speeds and perform a grid refinement study using a dG method of order 5, a Hermite method of order 9, and the Finite Volume method, all using the classical fourth order Runge-Kutta time stepping. For the Hermite method, we fix $(\max|u|)\Delta t /h = 0.3$, for the dG method, the time step is set as $\Delta t /h = 0.0625$ and for the Finite Volume method, the time step is set according to $(\max|u|)\Delta t/h = 0.9$. 

The $L_2$ norm of errors in the numerical solution $u_h$ are plotted against the different grid sizes for different methods, see Figure \ref{fig:APA}. In the stationary shock experiment, FV1 and FV2 use $(\alpha_{EV},\alpha_{\rm max}) = (0.7,0.5)$ and $(10,0.5)$ respectively, DG1 and DG2 use $(\alpha_{EV},\alpha_{\rm max}) = (1,0.25)$ and $(10,0.25)$ respectively, H1 and H2 use $(\alpha_{EV},\alpha_{\rm max}) = (1,0.4)$ and $(10,0.4)$ respectively, DGP and HP use $(s_0,\kappa,\epsilon_0) = ( -1, 2,0.5)$ and $(\log_{10}(1/256), 1,0.125)$ respectively.

The parameters for moving shock experiment are $(\alpha_{EV},\alpha_{\rm max}) = (0.7,0.5)$ and $(10,0.5)$ for FV1 and FV2 respectively, $(\alpha_{EV},\alpha_{\rm max}) = (1,0.25)$ and $ (10,0.25)$ for DG1 and DG2 respectively, $(\alpha_{EV},\alpha_{\rm max}) = (1,0.4)$ and $(10,0.4)$ for H1 and H2 respectively, $(s_0,\kappa,\epsilon_0) = (2\log_{10}(1/256) , 1,0.5)$ and $(\log_{10}(1/256), 1,0.125)$ for DGP and HP respectively.

To the left in Figure \ref{fig:APA}, we display convergence results for the stationary shock. In this case, the results indicate that all methods produce convergent solutions with roughly the same rates of convergence. The rate of convergence is limited by the smoothness of the solution but as can be seen in the same figure, the error levels are lower for the higher order methods. It is interesting to note that the smallest errors are observed for the computations using the smoothness-based sensor.

The results for the moving shock, displayed to the right in Figure \ref{fig:APA}, are quite different. Now, for the high order methods, we observe convergence only when we use the entropy viscosity with $\beta=1$. When we use the entropy viscosity with  $\beta=2$ or when we use the smoothness based sensor, the errors clearly saturate as the grid is refined. The errors for the low order Finite Volume method are still reduced with the grid size, independent of the scaling in the entropy viscosity method.

To understand why the convergence results obtained with $\beta=1$ and $\beta=2$ in the moving shock example do not agree, we study where the Lax-Friedrich viscosity $\nu_{\rm max}$ is activated in the vicinity of the shock. We know that when the viscosity is chosen to be just the Lax-Friedrich type viscosity, then under a suitable Courant number, the solution will converge to the correct vanishing viscosity solution of the conservation law  \cite{leveque2002finite}.

\begin{figure}[t]
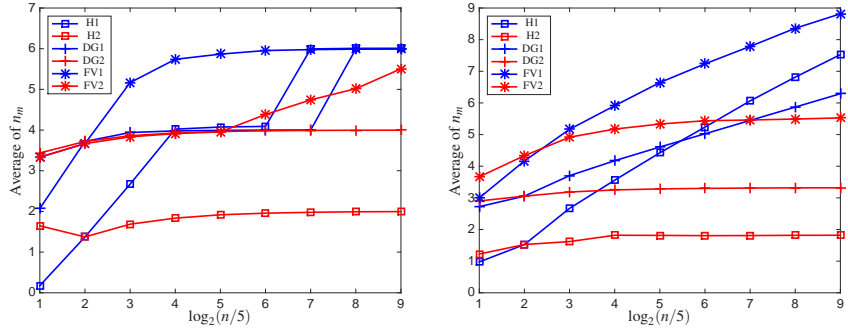

  \begin{center} 
  \psfrag{XXXXXXXX}[][][0.7][0]{$\log_2(n/5)$} 
  \psfrag{YYYYYYYY}[][][0.7][0]{Average of $n_m$}
  \psfrag{AAA1}[][][0.5][0]{H1} 
  \psfrag{AAA2}[][][0.5][0]{H2} 
  \psfrag{AAA3}[][][0.5][0]{DG1} 
  \psfrag{AAA4}[][][0.5][0]{DG2} 
  \psfrag{AAA5}[][][0.5][0]{FV1} 
  \psfrag{AAA6}[][][0.5][0]{FV2}
  \includegraphics[width=0.45\textwidth]{nogp_stationary_shock_1_2.eps}\hspace{0.5cm}
  \includegraphics[width=0.45\textwidth]{nogp_moveshock_1_2.eps}
  \caption{Average in time of $n_m$, number of elements using Lax-Friedrich viscosity $\nu_{\rm max}$, versus the number of elements ($n$). Left: stationary shock, right: moving shock. \label{fig:count_grid}}
  \end{center}
\end{figure}

It seems that the Lax-Friedrich viscosity is necessary in some neighborhood of the shock, and the size of this neighborhood becomes an important factor in the convergence of the solution to the moving shock problem. In Figure \ref{fig:count_grid}, we plot the average (in time) of the number of elements $n_m$ which use the Lax-Friedrich viscosity $\nu_{\rm max}$ as a function of total number of elements $n$ for the stationary shock (left) and for the moving shock (right). We see that $n_m$ is roughly constant for both $\beta=1$ and $\beta=2$ in the stationary shock. In the moving shock problem, $n_m$ stays constant for $\beta=2$ as in the stationary shock, but grows slowly for $\beta=1$ (note the log-scale). While the growth in $n_m$ is irrelevant in the convergence in the stationary shock example, it seems to play an important role in determining the convergence in the moving shock example.

\subsection{Sinusoidal to N wave}
Next, we consider the  smooth 2-periodic initial data 
\be
u(x,0) = -\sin(\pi x)+0.5,
\ee
which develops into a single N wave.

\begin{figure}[t]
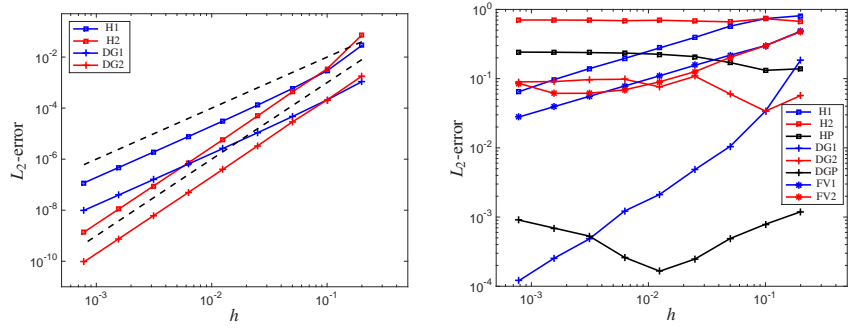

  \begin{center} 
  \psfrag{XXXXXXXX}[][][0.8][0]{$h$} 
  \psfrag{YYYYYYYY}[][][0.8][0]{$L_2$-error} 
  \psfrag{AAA1}[][][0.5][0]{H1} 
  \psfrag{AAA2}[][][0.5][0]{H2} 
  \psfrag{AAA3}[][][0.5][0]{HP} 
  \psfrag{AAA4}[][][0.5][0]{DG1} 
  \psfrag{AAA5}[][][0.5][0]{DG2} 
  \psfrag{AAA6}[][][0.5][0]{DGP} 
  \psfrag{AAA7}[][][0.5][0]{FV1} 
  \psfrag{AAA8}[][][0.5][0]{FV2} 
  \includegraphics[width=0.45\textwidth]{smooth_sine_wave.eps}\hspace{0.5cm}
  \includegraphics[width=0.45\textwidth]{moving_sine.eps}
  \caption{Convergence of the different methods for a smooth initial data, left: before shock forms, right: after shock forms. The dashed lines are $h^2$ and $h^3$. \label{fig:smooth_sine}}
  \end{center}
\end{figure}

In Figure \ref{fig:smooth_sine}, we present the $L_2$ norm of the errors at $t=0.1$ before the shock forms (left) and at $t=1$ after the shock forms (right). The spatial and temporal discretization of the PDE itself is performed with a high order method, so rate of convergence that we observe in Figure \ref{fig:smooth_sine} is limited by either the discretization of the artificial viscosity or the smoothness of the solution, whichever is more restrictive.

For this N-wave experiment, $FV1$ and $FV2$ use ($\alpha_{EV},\alpha_{\rm max}) = (2,0.5)$ and $(20,0.5)$ respectively, $DG1$ and $DG2$ use ($\alpha_{EV},\alpha_{\rm max}) = (0.1,0.125)$ and $(1,0.125)$ respectively, $H1$ and $H2$ use ($\alpha_{EV},\alpha_{\rm max}) = (0.4,0.4)$ and $(5,0.4)$, $DGP$ and $HP$ use $(s_0,\kappa,\epsilon_0) = (2\log_{10}(1/256), 2,0.05)$ and $(\log_{10}(1/256), 1,0.125)$ respectively.

The discretization of the entropy residual $r_{EV}$ is only first order due to the use of backward-Euler, so we expect the entropy-based viscosity $\nu_{EV}$ to be $(\beta+1)^{th}$ accurate, i.e. $2^{nd}$ order when $\beta=1$ or $3^{rd}$ order when $\beta=2$. This analysis agrees with the convergence plot to the left in Figure \ref{fig:smooth_sine}. To the right, we observe the same phenomena as in the moving shock example described in Section \ref{sec:single_shock}. We also note that the shock in this sinusoidal wave is also moving.

%
%

\subsection{Shocks of different size} \label{sec:two_shocks}

To complement the analysis in Section \ref{sec:viscosity_analysis}, we next consider a problem with a big shock and a small shock on the same simulation. According to the analysis, the entropy viscosity will capture the small shock when $\beta=1$, but not when $\beta=2$. In this setup, we start with an existing shock of size $\Delta u_1=0.5$ and a small sinusoidal wave that develops into an N-wave of size $\Delta u_2=0.2$. Thus, we consider Burgers' equation on $[-1,5]$ with initial data 
\be
u(x,0) =
\begin{cases}
  \ \ \; 0          & x \in [-1,-0.5),\\
 -0.1 \sin (2\pi x) & x \in [-0.5,0.5), \\
  \ \ \; 0          & x \in [0.5,4.5), \\
  -0.5         & x \in [4.5,5],
\end{cases}
\ee
and fixed boundary condition $u(-1,t) = 0$ and $u(5,t) = -0.5$.

The solution initially consists of a shock and a smooth sine wave, which are placed far away from each other so they never interact. Over time, the sinusoidal wave develops into a N-wave. In Figure \ref{fig:move_pert2}, we present the numerical solutions at time $t = 2$ for different grid resolutions, obtained with a Hermite method of order 9 and dG method of order 5. In these plots, we can see that the shock is resolved for both values of $\beta$, however, the N-wave comes with some overshoots when $\beta=2$ for all the finer grid resolutions, see Figure \ref{fig:move_pert2}.

For this two-shock experiment, $DG1$ and $DG2$ use ($\alpha_{EV},\alpha_{\rm max}) = (0.5,0.25)$ and $(10,0.25)$ respectively, $H1$ and $H2$ use ($\alpha_{EV},\alpha_{\rm max}) = (1,0.125)$ and $(50,0.125)$ respectively.

\begin{figure}[htb]
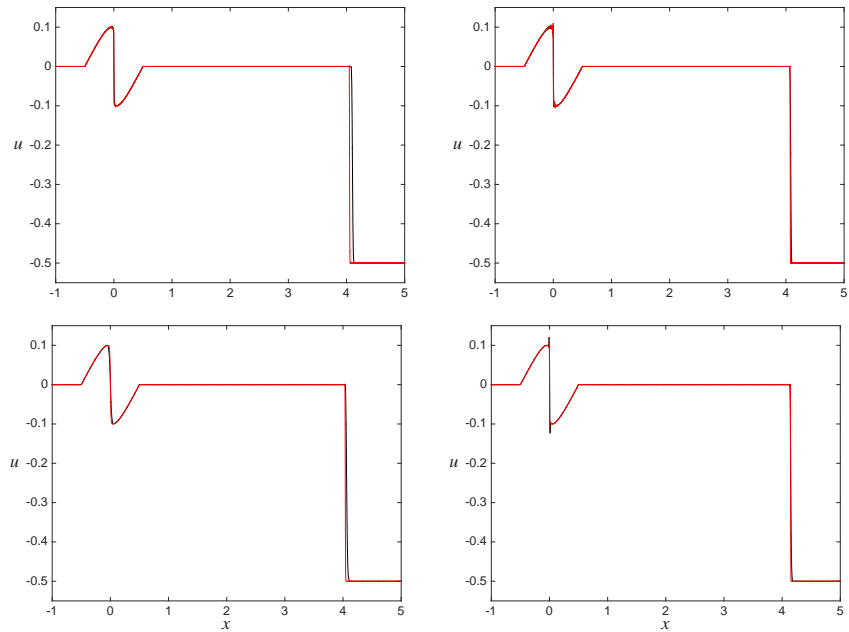

  \begin{center} 
  \psfrag{XXXX}[][][0.8][0]{$x$} 
  \psfrag{YYYY}[][][0.8][270]{$u$} 
  \includegraphics[width=0.45\textwidth]{pert_full_dg_b1.eps}\hspace{0.5cm}
  \includegraphics[width=0.45\textwidth]{pert_full_dg_b2.eps}\hspace{0.1cm} \\ \vspace{0.1cm}
  \includegraphics[width=0.45\textwidth]{pert_full_H_b1.eps}\hspace{0.5cm}
  \includegraphics[width=0.45\textwidth]{pert_full_H_b2.eps}\hspace{0.1cm}
  \caption{Effect of the choice of scaling on a small perturbation near a larger shock. The results in the left and right column are for $\beta = 1$ and $\beta =2$ respectively. The upper figures display the results for the dG method and the lower figures display the results for the Hermite method. The black curve is for a simulation using 320 elements and the black uses 2560.  \label{fig:move_pert}}
  \end{center}
\end{figure}

\begin{figure}[htb]
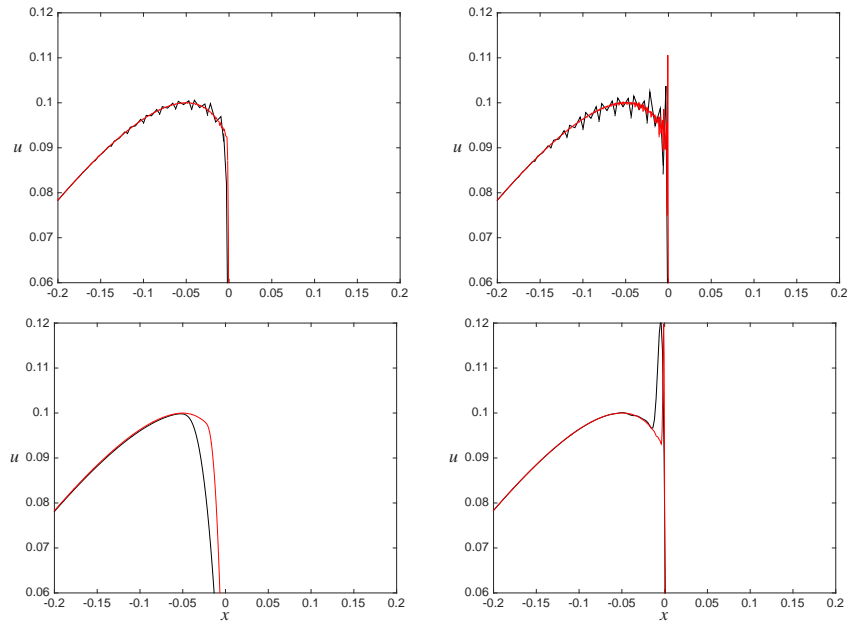

  \begin{center} 
  \psfrag{XXXX}[][][0.8][0]{$x$} 
  \psfrag{YYYY}[][][0.8][270]{$u$} 
  \includegraphics[width=0.45\textwidth]{pert_zoom2_dg_b1.eps}\hspace{0.5cm}
  \includegraphics[width=0.45\textwidth]{pert_zoom2_dg_b2.eps}\hspace{0.1cm} \\
  \includegraphics[width=0.45\textwidth]{pert_zoom2_H_b1.eps}\hspace{0.5cm}
  \includegraphics[width=0.45\textwidth]{pert_zoom2_H_b2.eps}\hspace{0.1cm}
  \caption{Effect of the choice of scaling on a small perturbation near a larger shock. Same as in Figure \ref{fig:move_pert} but zoomed in. \label{fig:move_pert2}}
  \end{center}
\end{figure}

Because the magnitude of this N-wave is small, the entropy residual at the N-wave is relatively small compared to that at the existing shock. On one hand, $\beta=1$ results are free from overshoots when the grid is refined, but $\beta=2$ results do have overshoots, see Figure \ref{fig:move_pert}-\ref{fig:move_pert2}.


\section{Conclusion}
In summary, we have performed a convergence study for Burgers' equation with various initial data. We demonstrated that the entropy viscosity method with $\beta=2$ does not produce convergent results (fixing the parameters $\alpha_{EV}$ and $\alpha_{\rm max}$) in the cases where the shock is moving or more than one shock is present. Therefore, we recommend readers to use $\beta=1$; to achieve desired accuracy or better rate of convergence, use a higher order approximation of the residual.

\bibliographystyle{siam}
\bibliography{kornelus}

\end{document}